\documentclass[fleqn]{mat01}
\usepackage{times,mathtimy,amssymb,latexsym,amsbsy,amsbsy,epsf,epsfig}%rangecite,
\begin{document}

\setcounter{page}{251}
\firstpage{251}

\newtheorem{theore}{Theorem}
\renewcommand\thetheore{\arabic{section}.\arabic{theore}}
\newtheorem{theor}[theore]{\bf Theorem}
\newtheorem{rem}[theore]{Remark}
\newtheorem{propo}[theore]{\rm PROPOSITION}
\newtheorem{lem}[theore]{Lemma}
\newtheorem{definit}[theore]{\rm DEFINITION}
\newtheorem{coro}[theore]{\rm COROLLARY}
\newtheorem{exampl}[theore]{Example}
\newtheorem{case}{Case}

\def\corol{\trivlist \item[\hskip \labelsep{COROLLARY.}]}
\def\noteproof{\trivlist \item[\hskip \labelsep{\it Note added in Proof.}]}

\renewcommand{\theequation}{\thesection\arabic{equation}}

\newcommand{\Z}{\mathbb{Z}}
\newcommand{\R}{\mathbb{R}}
\newcommand{\C}{\mathbb{C}}
\newcommand{\N}{\mathcal{N}}
\newcommand{\M}{\mathcal{M}}
\newcommand{\K}{\mathcal{K}}
\newcommand{\Zp}{\mathbb{Z}/p\mathbb{Z}}
\renewcommand{\S}{\mathcal{S}}
\renewcommand{\P}{\mathcal{P}}
\renewcommand{\O}{\mathcal{O}}
\renewcommand{\H}{{\mathbb{H}}}
\newcommand{\F}{\mathcal{F}}
\newcommand{\Zz}{\mathbb{Z}/2\mathbb{Z}}
\newcommand{\Zq}{\mathbb{Z}/q\mathbb{Z}}
\newcommand{\Zn}{\mathbb{Z}/n\mathbb{Z}}
\newcommand{\sm}{\setminus}
\newcommand{\tr}{\triangle}
\renewcommand{\i}{\bm{i}}
\renewcommand{\j}{\bm{j}}
\renewcommand{\k}{\bm{k}}
\newcommand{\del}{\partial}
\newcommand{\tF}{\tilde F}
\newcommand{\tC}{\tilde C}
\newcommand{\tS}{\tilde\Sigma}
\newcommand{\nbd}{\mathcal{N}}

\makeatletter
\newcommand\abb{\@setfontsize\abb{7.5}{9}}
\def\artpath#1{\def\@artpath{#1}}
\makeatother \artpath{C:/mathsci-arxiv/aug2005/}

\title{Homeomorphisms and the homology of non-orientable surfaces}

\markboth{Siddhartha Gadgil and Dishant Pancholi}{Homeomorphisms and the homology of
non-orientable surfaces}

\author{SIDDHARTHA GADGIL and DISHANT PANCHOLI$^{*}$}

\address{Stat Math Unit, Indian Statistical Institute,
Bangalore~560~059, India\\
\noindent $^{*}$School of Mathematics, Tata Institute of
Fundamental Research,
Mumbai~400~005, India\\
\noindent E-mail: gadgil@isibang.ac.in; dishant@math.tifr.res.in}

\volume{115}

\mon{August}

\parts{3}

\pubyear{2005}

\Date{MS received 8 February 2005; revised 3 May 2005}

\begin{abstract}
We show that, for a closed non-orientable surface $F$, an
automorphism of $H_1(F,\Z)$ is induced by a homeomorphism of $F$
if and only if it preserves the (mod $2$) intersection pairing. We
shall also prove the corresponding result for punctured surfaces.
\end{abstract}


\maketitle

\section{Introduction}

Let $F$ be a closed, non-orientable surface. A homeomorphism
$f\hbox{\rm :}\ F\rightarrow  F$ induces an automorphism on
homology $f_*\hbox{\rm :}\ H_1(F,\Z)\rightarrow  H_1(F,\Z)$.
Further, any automorphism $\varphi\hbox{\rm :}\
H_1(F,\Z)\rightarrow H_1(F,\Z)$ in turn induces an automorphism
with $\Zz$-coefficients $\bar\varphi\hbox{\rm :}\
H_1(F,\Zz)\rightarrow H_1(F,\Zz)$. If $\varphi=f_*$ for a
homeomorphism $f$, then $\bar\varphi$ also preserves the (mod~$2$)
intersection pairing on homology.

Our main result is that, for an automorphism $\varphi\hbox{\rm :}\
H_1(F,\Z)\rightarrow  H_1(F,\Z)$, if the induced automorphism
$\bar\varphi\hbox{\rm :}\ H_1(F,\Zz)\rightarrow  H_1(F,\Zz)$
preserves the (mod~$2$) intersection pairing, then $\varphi$ is
induced by a homeomorphism of $F$.

\begin{theor}[\!]\label{main}
Let $\varphi\hbox{\rm :}\ H_1(F,\Z)\rightarrow H_1(F,\Z)$ be an
automorphism. If the induced automorphism $\bar\varphi\hbox{\rm
:}\ H_1(F,\Zz)\rightarrow H_1(F,\Zz)$ preserves the {\rm
(}\hbox{\rm mod 2}{\rm )} intersection pairing{\rm ,} then
$\varphi$ is induced by a homeomorphism of $F$.
\end{theor}

We have a natural homomorphism $\hbox{Aut}(H_1(F,\Z))\rightarrow
\hbox{Aut}(H_1(F,\Zz))$. Let $\K$ denote the kernel of this
homomorphism, so that we have an exact sequence
\begin{equation*}
1\rightarrow \K\rightarrow \hbox{Aut}(H_1(F,\Z))\rightarrow
\hbox{Aut}(H_1(F,\Zz))\rightarrow 1.
\end{equation*}
Observe that elements of $\K$ automatically preserve the
intersection pairing. We shall show that every element of $\K$ is
induced by a homeomorphism of $F$. Further, we shall show that an
element of $\hbox{Aut}(H_1(F,\Zz))$ is induced by a homeomorphism
of $F$ if and only if it preserves the intersection pairing.
Theorem~\ref{main} follows immediately from these results.

\begin{theor}[\!]\label{ker}
Suppose $\varphi\hbox{\rm :}\ H_1(F,\Z)\rightarrow H_1(F,\Z)$ is
an automorphism which induces the identity on $H_1(F,\Zz)$. Then
$\varphi$ is induced by a homeomorphism of $F$.
\end{theor}

\begin{theor}[\!]\label{mod2}
Let $F_1$ and $F_2$ be closed{\rm ,} non-orientable surfaces.
Suppose that\break $\psi\hbox{\rm :}\ \ H_1(F_1,\Zz)\rightarrow
H_1(F_2,\Zz)$ is an isomorphism which preserves the intersection
pairing. Then $\psi$ is induced by a homeomorphism $f\hbox{\rm :}\
F_1\rightarrow F_2$.
\end{theor}

We also consider the case of a compact non-orientable surface $F$
with boundary. In this case an automorphism of $H_1(F,\Z)$ induced
by a homeomorphism of $F$ permutes (up to sign) the elements
representing the boundary components. We shall show that all
automorphisms of $H_1(F,\Z)$ which satisfy this additional
condition are induced by homeomorphisms. Other results regarding the
homeomorphisms of non-orientable surfaces have been obtained by many
authors, for instance \cite{BC,Ko,Li}.
\vspace{-.3pc}

\section{Preliminaries}

Let $F$ be a closed, non-orientable surface with $\chi(F)=2-n$ and
let $\hat F$ be obtained from $F$ by deleting the interior of a
disc. Then $F$ is the connected sum of $n$ projective planes
$\P_i$ and $\hat F$ is the boundary-connected sum of $n$
corresponding M\"obius bands $\M_i$. Let $\gamma_i$ denote the
central circle of $\M_i$ and let $\alpha_i=[\gamma_i]\in H_1(\hat
F,\Z)$ be the corresponding elements in homology. Then $H_1(\hat
F,\Z)\ \cong \Z^n$ with basis $\alpha_i$ and $H_1(F,\Z)$
is the quotient $H_1(\hat F,\Z)/{\langle
2\Sigma_i\alpha_i\rangle}$.

We shall need the following elementary algebraic lemma.

\setcounter{theore}{0}
\begin{lem}\label{lift}
Any automorphism $\varphi\hbox{\rm :}\ H_1(F,\Z)\rightarrow
H_1(F,\Z)$ lifts to an automorphism $\tilde\varphi\hbox{\rm :}\
H_1(\hat F,\Z)\rightarrow  H_1(\hat F,\Z)$ such that
$\tilde\varphi(\sum_i\alpha_i)=\sum_i\alpha_i$.
\end{lem}

\begin{proof}
Consider the basis of $H_1(\hat F,\Z)$ given by
$e_1=\alpha_1,\ldots, e_{n-1}=\alpha_{n-1}$,
$e_n=\alpha_1+\cdots+\alpha_n$ and let $[e_j]$ be the
corresponding generators of $H_1(F,\Z)$. Observe that $[e_n]$ is
the unique element of order $2$ in $H_1(F,\Z)$, and hence
$\varphi([e_n])=[e_n]$. Thus, we can define
$\tilde\varphi(e_n)=e_n$. For $1\leq j\leq n-1$, pick an arbitrary
lift $h_j$ of $\varphi(e_j)$ and set $\tilde\varphi(e_j)=h_j$.

Observe that $H_1(\hat F,\Z)/{\langle e_n\rangle}\ \cong
H_1(F,\Z)/{\langle [e_n]\rangle}$. Further, as
$\tilde\varphi(e_n)=e_n$ we have an induced map on $H_1(\hat
F,\Z)/{\langle e_n\rangle}$ which agrees with the quotient map
induced by $\varphi$ on $H_1(F,\Z)/{\langle [e_n]\rangle}$ (which
exists as $\varphi([e_n])=[e_n]$) under the natural identification
of these groups. As $\varphi$ is an isomorphism, so is the induced
quotient map on $H_1(F,\Z)/{\langle [e_n]\rangle}$, and hence the
map induced by $\tilde\varphi$ on $H_1(\hat F,\Z)/{\langle
e_n\rangle}$.

Thus, $\tilde\varphi$ induces an isomorphism on the quotient
$H_1(\hat F,\Z)/{\langle e_n\rangle}$ as well as the kernel
$\langle e_n\rangle$ of the quotient map. By the five lemma,
$\tilde\varphi$ is an isomorphism.\hfill $\Box$
\end{proof}

Henceforth, given an automorphism $\varphi$ as above, we shall
assume that a lift has been chosen as in the lemma. Observe that a
homeomorphism of $\hat F$ induces a homeomorphism of $F$. Hence it
suffices to construct a homeomorphism of $\hat F$ inducing
$\tilde\varphi$. Note that the intersection pairing is preserved
by $\tilde\varphi$ as it only depends on the induced map on
homology with $\Zz$-coefficients.\vspace{-.3pc}

\section{Automorphisms in $\pmb{\K}$}

In this section we prove Theorem~\ref{ker}. Let $\varphi\hbox{\rm
:}\ H_1(F,\Z)\rightarrow  H_1(F,\Z)$ be as in the hypothesis. As
in Lemma~\ref{lift}, we can lift $\varphi$ to an automorphism of
$H_1(\hat F,\Z)$ fixing $\sum_i\alpha_i$. We shall denote this
lift also by $\varphi$. We shall construct a homeomorphism of
$\hat F$ inducing this automorphism.

Our strategy is to use {\it elementary automorphisms} $e_{ij}$,
$1\leq i,j\leq n$, which are induced by homeomorphisms $g_{ij}$.
Observe that, for $1\leq i,j\leq n$, the automorphism $\varphi$ is
induced by a homeomorphism if and only if $e_{ij}\circ \varphi$ is
induced by a homeomorphism (as $e_{ij}$ is induced by a
homeomorphism). Thus we can replace $\varphi$ with
$e_{ij}\circ\varphi$. We call this an {\it elementary move}. For
$\varphi$ preserving the intersection pairing, we shall find a
sequence of elementary moves such that on performing these moves
we obtain the identity automorphism, which is obviously induced by
a homeomorphism (namely the identity). This will prove the result.

\begin{figure}[t]
%\begin{center}
%\parbox{0.55\textwidth}{\epsfxsize=0.55\textwidth
\centerline{\epsfxsize=9cm\epsfbox{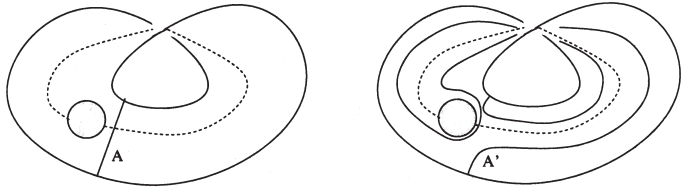}}\vspace{.8pc}
\caption{Cross-cap slide.}\label{cross}\vspace{1pc}
%\end{center}
\end{figure}

%\begin{figure}%1
%\hskip 4pc {\epsfxsize=10cm\epsfbox{pm2508-1.eps}}\vspace{-.5pc}
%\caption{Cross-cap slide.}\label{cross}\vspace{1pc}
%\end{figure}

\setcounter{theore}{0}
\begin{lem}
There are homeomorphisms $g_{ij}$ of $\hat F$ so that if $e_{ij}$
is the induced automorphism on $H_1(\hat F,\Z)${\rm ,} then
$e_{ij}(\alpha_i)=\alpha_i+2\alpha_j${\rm ,}
$e_{ij}(\alpha_j)=-\alpha_j$ and $e_{ij}(\alpha_k)=\alpha_k$ for
$k\neq i,j$.
\end{lem}

\begin{proof}
We shall use {\it cross-cap slides}~\cite{Li,Li1} of the surface
$F$. Namely, suppose $\alpha$ is an orientation reversing simple
closed curve on a surface $S'$ and $D$ is a small disc centered
around a point on $\alpha$. Let $S$ be the surface obtained by
replacing $D$ by a M\"obius band. Consider a homeomorphism $f'$ of
$S'$ which is the identity outside a neighbourhood of $\alpha$ and
which is obtained by dragging $D$ once around $\alpha$ so that $D$
is mapped to itself. By construction this extends to a
homeomorphism $f$ of $S$, which we call a cross-cap slide. In
figure~\ref{cross}, the arc $A$ in the M\"obius band $\M$ on the
left-hand side is mapped to the arc $A'$ in the M\"obius band
$\M'$ on the right-hand side and the homeomorphism is the identity
in a neighbourhood of the boundary.

We define $g_{ij}$ as the {\it cross-cap slide} of $\M_{j}$ around
the curve $\gamma_{i}$. Note that the M\"obius band $\M_{j}$ is
mapped to itself, but, as $\gamma_{i}$ is orientation reversing,
the map on the M\"obius band takes $\alpha_j$ to $-\alpha_j$.
Further for any $k$ different from $i$ and $j$, the cross-cap
slide fixes $\gamma_j$, hence $\alpha_j$. Finally, in
figure~\ref{cross} (where we regard $\M$ as a neighbourhood of
$\alpha_j$), if $B$ is a curve in the boundary of $\M$ joining the
endpoints of $A$, then $[A\cup B]=\alpha_i$ and $[A'\cup
B]=g_{ij}(\alpha_i)$. It is easy to see that $[A\cup B]-[A'\cup
A]=[A\cup A']$ is homologous to the boundary $2\alpha_j$ of the
cross-cap. Thus, $e_{ij}(\alpha_i)=\alpha_i+2\alpha_j$.\hfill
$\Box$
\end{proof}

\begin{lem}
There exists a sequence of elementary moves $e_{ij}$ taking
$\varphi$ to the identity.
\end{lem}

\begin{proof}
Let $\varphi$ be represented by a matrix $A=(a_{ij})$ with respect
to the basis $\alpha_i$. Then $A\equiv I(\hbox{mod}\ 2)$. As
$\varphi$ fixes $\sum_i\alpha_i$, for every $i, \sum_j a_{ij}=1$.
Observe that on performing the elementary move $e_{ij}$, the $i$th
column $A_{*i}$ of $A$ is replaced by $A_{*i}+2A_{*j}$, the $j$th
column is replaced by $-A_{*j}$ and the other columns of $A$ are
unchanged.

We first use the elementary moves $e_{ij}$ to reduce the first row
$A_{1*}$ to $[1,0, 0, \ldots ,0]$. To do this, we define a
complexity $C_{1}(A)$ of $A$ as $|a_{11}| + |a_{12}| + \cdots +
|a_{1n}|$.

Observe that if $a_{1k}$ and $a_{1l}$ are both non-zero, have
different signs and $|a_{1k}|>|a_{1l}|$, $e_{kl}$ reduces the
complexity $C_1(A)$. As $a_{11}$ is odd and $a_{1j}$ is even for
$j>1$, we know that $a_{11}\neq a_{1j}$ for every $j>1$. Further,
as $\sum_j a_{1j}=1$, unless $a_{11}$ is $1$ and $a_{1j}=0$ for $j
\neq 1$, there exists a $j>1$ such that $a_{11}$ and $a_{1j}$ are
of opposite signs (and both non-zero). Thus we can reduce
complexity by performing an elementary operation. By iterating
this finitely many times, we reduce the first row to $[1,0,
0,\ldots ,0]$.

Next, suppose $i>1$ and the rows $A_{1*}$, $A_{2*}$,\ldots,
$A_{(i-1)*}$ are the unit vectors $e_1$, $e_2$,\ldots,
$e_{(i-1)*}$. We shall transform the $i$th row to
$[0,0,\ldots,1,0,\ldots,0]$ without changing the earlier rows.

First we shall transform the row $A_{i*}$ to a row of the form
$[*,*,\ldots,1,0,\ldots,0]$ (i.e., with the first $i-1$ entries
arbitrary) by performing elementary moves $e_{ij}$. To do this, we
define a complexity $C_{i}(A) = \sum_j |a_{ij}|, j\geq i$.

Observe that, for $k \geq i$, the elementary operation $e_{1k}$
changes the sign of $a_{ik}$, does not alter $a_{im}$ for $m \neq
k, m \geq i$ and does not change first $i-1$ rows.  By such
operations we can ensure that $a_{ii}>0$ and $a_{ij}\leq 0$ for
$j>i$ without changing the complexity.

As before, $a_{ii} \neq a_{ij}$ for $j>i$ (as $a_{ii}$ is odd and
$a_{ij}$ is even) and (using operations $e_{1k}$ if necessary)
$a_{ii}$ and $a_{ij}$ have different signs. Hence, unless $a_{ij}
= 0$ for $j>i$ we can reduce the complexity using either $e_{ij}$
or $e_{ji}$, without altering the first $i$ rows. Thus we can
reduce $A_{i*}$ to a vector of the form
$[*,\ldots,*,m,0,\ldots,0]$.

Now $A$ is a block lower triangular matrix with $a_{ii}$ as a
diagonal entry. As $A$ is invertible it follows that $m=a_{ii}=\pm
1$.

We define another complexity $C'_i(A)=\sum_{i\geq j} |a_{ij}|$. As
$\sum_j a_{ij}=1$ and $a_{ii}=\pm 1$, unless $A_{i*}$ is a unit
vector we can find as before an operation $e_{ji}, j<i$, which
reduces this complexity (without changing the first $(i-1)$ rows).
Hence after finitely many steps the $i$th row is reduced to a unit
vector. By applying these moves for $i=2,3,\ldots, n$, we are
done.\hfill $\Box$
\end{proof}

\section{Automorphisms of $\pmb{H_1(F,\Zz)}$}

We now prove Theorem~\ref{mod2}. We shall proceed by induction on
$n$. In the case when $n=1$ the result is obvious. We henceforth
assume that $n$ is greater than $1$.

We first make some observations. For a surface $S$, any element
$\alpha$ of $H_{1}(S,\Zz)$ can be represented by a simple closed
curve. The curve $\alpha$ is orientation reversing if and only if
$\alpha\cdot\alpha = 1$. The surface is non-orientable if and only
if there exist $\alpha \in H_{1}(S,\Zz)$ with $\alpha \cdot\alpha =1$.

As before, let $F_{1}$ be the connected sum of
$\P_{1},\P_{2},\ldots,\P_{n}$, where $\P_i$ denotes a projective
plane and $\M_{i}$ denotes the corresponding M\"obius band. Let
$\alpha_{1},\alpha_{2},\ldots,\alpha_{n}$ and
$\gamma_1$,\ldots,$\gamma_n$ be as before.

Let $\psi$ be as in the hypothesis. Let $\beta_i=\psi(\alpha_{i})$
and let $C$ be a simple close curve that represents $\beta_{1}$.
As $\beta_{1}\cdot\beta_{1} = \alpha_{1}\cdot\alpha_{1}=1$, $C$ is
orientation reversing (as is $\gamma_{1}$). Hence regular
neighbourhoods of $C$ and $\gamma_{1}$ are M\"obius bands.

Let $\hat {F'_{1}} = F_{1} - \hbox{int}(\N(\gamma_{1}))$ and $\hat
F'_{2}=F_{2} - \hbox{int}(\N(C))$. Let $ F'_{1} = \hat F'_{1}
\cup D^{2}$ and $F'_{2} = \hat F'_{2} \cup D^{2}$ be closed
surfaces obtained by capping off $\hat F_i$.

Observe that the surface $F_{1}'$ is non-orientable as $n \geq 2$
and $\gamma_{2}$ is an orientation reversing curve on it. Now
since $\psi$ preserves the intersection pairing it takes
orthonormal basis of $H_{1}(F_{1},\Zz)$ to orthonormal basis of
$H_{1}(F_{2},\Zz)$. It follows that $\beta_{j}\, \cdot\, \beta_{j} =1$
for every $j$. Further, by a Mayer--Vietoris argument,
$H_1(F_i,\Zz)=\Zz\,\oplus\, H_1(F_i',\Zz)$, with the decomposition
being orthogonal and the component $\Zz$ in $H_1(F_1,\Zz)$
(respectively $H_1(F_2,\Zz)$) is spanned by $\alpha_1$
(respectively $\beta_1$). As $\psi$ preserves the intersection
pairing, it follows that $\psi$ induces an isomorphism
$\psi\hbox{\rm :}\ H_1(F_1,\Zz)\rightarrow  H_1(F_2,\Zz)$.

Hence if $C_{2}$ is a curve in $F'_{2}$ representing $\beta_{2}$
in $H_{1}(F_{2},\Zz)$, then $C_{2}$ is orientation reversing and
hence $F'_2$ is non-orientable. Also, we have seen that the map
$\psi$ induces an isomorphism from $H_{1}(F_{1}',\Zz)$ to
$H_{1}(F_{2}',\Zz)$.  By the induction hypothesis such a map is
induced by a homoeomorphism $f'\hbox{\rm :}\ F'_{1} \rightarrow
F'_{2}$.

Note that $F_1$ (respectively $F_2$) is obtained from $F'_1$
(respectively $F'_2$) by deleting the interior of a disc $D_1$
(respectively $D_2$) and gluing in $\N(\gamma_1)$ (respectively
$\N(C)$). We can modify $f'$ so that $f'(D_1)=D_2$. On
$F_1-\hbox{int}(D_1)$ we define $f=f'$. This restricts to a
homeomorphism mapping $\del \N(\gamma_1)$ to $\del \N(C)$, which
extends to a homeomorphism mapping $\N(\gamma_1)$ to $\N(C)$. As
$f|_{\N(\gamma_1)}\hbox{\rm :}\ \N(\gamma_1)\rightarrow  \N(C)$ is
a homeomorphism, it maps the generator $\alpha_1$ of
$H_1(\N(\gamma_1),\Z)= \Z$ to a generator $\pm\beta$ of
$H_1(\N(C),\Z)= \Z$. Thus with mod~2 coefficients, $f_*=\varphi$
as required.

\section{An algebraic corollary}

We shall deduce from Theorem~\ref{main} and a theorem of
Lickorish~\cite{Li} a purely algebraic corollary. While this has a
straightforward algebraic proof (and is presumably well-known), it
may still be of interest to see its relation to topology.

Let $V=(\Zz)^n$ be a vector space over $\Zz$ and let $\{e_j\}$ be
the standard basis of $V$. Consider the standard inner product
$\langle (x_i),(y_i)\rangle=\sum_{i} x_i y_i$. Let $\O$ be the
group of automorphisms of $V$ that preserve the inner product. We
shall show that $\O$ is generated by certain involutions.

Namely, let $1\leq i_1<i_2<\cdots<i_{2k}\leq n$ be $2k$ integers
between $1$ and $n$. We define an element $R=R(i_1,\ldots,i_{2k})$
to be the transformation defined by
\begin{align*}
R(e_{i_j}) &= \sum_{l\neq j}e_{i_l},\\
R(e_j) &= e_j, j\neq i_1,i_2,\ldots i_{2k}.
\end{align*}

\setcounter{theore}{0}
\begin{theor}[\!]\label{alg}
The group $\O$ is generated by the involutions
$R(i_1,\ldots,i_{2k})$.
\end{theor}

\begin{proof}
We identify $V$ with $H_1(F,\Zz)$ for a non-orientable surface $F$
and identify the basis elements $e_i$ with $\alpha_i$. Under this
identification, the bilinear pairing on $V$ corresponds to the
intersection pairing. We shall see that the transformations
$R(i_1,\ldots,i_{2k})$ correspond to the action of Dehn twists on
$H_1(F,\Zz)$, where we identify the generators $e_i$ with
$\alpha_i$. First note that any element $\gamma$ of $H_1(F,\Zz)$
can be expressed as $\gamma=\alpha_{i_1}+\cdots+\alpha_{i_m}$.
Observe that a simple closed curve $C$ representing $\gamma$ is
orientation preserving if and only if $\gamma\cdot\gamma=0$, which
is equivalent to $m$ being even.

Now let $C$ be an orientation preserving curve on $F$ and consider
the Dehn twist $\tau$ about $C$. Let $\gamma=[C]\in H_1(F,\Zz)$ be
the element represented by $C$. By the above (as $\gamma = \alpha_{i_1}
+ \dots + \alpha_{i_m}$ and $m$ is even), we can express
$\gamma$ as $\gamma=\alpha_{i_1}+\cdots+\alpha_{i_{2k}}$. If
$\alpha$ is another element of $H_1(F,\Zz)$ and
$\alpha\cdot\gamma$ is the (mod~2) intersection number, then (with
mod~2 coefficients) $\tau_*(\alpha)=\alpha+\gamma$. It is easy to
see that $\tau_*=R(i_1,\ldots,i_{2k})$. Note that
$\tau^2_*(\alpha)=\alpha+2\gamma=\alpha$, hence
$\tau_*=R(i_1,\ldots,i_{2k})$ is an involution as claimed.

Now, by Theorem~\ref{mod2}, any element $\phi\in\O$ is induced by
a homeomorphism $f$ of $F$. Further, by a theorem of
Lickorish~\cite{Li}, $f$ is homotopic to a composition of Dehn
twists and cross-cap slides. We have seen that Dehn twists induce
the automorphisms $\tau_*=R(i_1,\ldots,i_{2k})$ on $V$. It is easy
to see that cross-cap slides induce the identity on $H_1(F,\Zz)$.
Thus $\phi$ is a composition of elements of the form
$\tau_*=R(i_1,\ldots,i_{2k})$ as claimed.\hfill $\Box$
\end{proof}

\begin{rem}
{\rm We can alternatively deduce Theorem~\ref{mod2} from
Theorem~\ref{alg} as the generators of $\O$ can be represented by
homeomorphisms (namely Dehn twists).}\vspace{-.5pc}
\end{rem}

\section{Punctured surfaces}

Let $F$ be a compact non-orientable surface with $m$ boundary
components and let $\beta_j\in H_1(F,\Z)$, $1\leq j\leq m$, be
elements representing the boundary curves. A homeomorphism
$f\hbox{\rm :}\ F\rightarrow  F$ induces an automorphism
$\varphi=f_*$ of $H_1(F,\Z)$. Furthermore, as boundary components
of $F$ are mapped to boundary components by $f$ (possibly
reversing orientations), for some permutation $\sigma$ of
$\{1,\ldots,m\}$ and some constants $\epsilon_j=\pm 1$,
$\varphi(\beta_j)=\epsilon_j\beta_{\sigma(j)}$, for all $j$,
$1\leq j\leq m$.

We show that conversely any automorphism $\varphi$ that preserves
the (mod~2) intersection pairing and takes boundary components to
boundary components is induced by a homeomorphism.

\setcounter{theore}{0}
\begin{theor}[\!]\label{punct}
Let $F$ be a compact non-orientable surface with $m$ boundary
components and let $\varphi$ be an automorphism of $H_1(F,\Zz)$
that preserves the {\rm (}mod~2{\rm )} intersection pairing.
Suppose for some permutation $\sigma$ of $\{1,\ldots,m\}$ and some
constants $\epsilon_j=\pm 1${\rm ,} we have
$\varphi(\beta_j)=\epsilon_j\beta_{\sigma(j)}${\rm ,} for all
$1\leq j\leq m$. Then $\varphi$ is induced by a homeomorphism of
$F$.
\end{theor}

\begin{proof}
Let $\bar F$ be obtained from $F$ by attaching discs to all the
boundary components. Then we can assume that $F$ has been obtained
from $\bar F$ by deleting the interiors of $m$ discs $D_1$,\ldots
$D_m$, all of which are contained in a disc $E\subset\bar F$.
Further we can assume that the central curves $\gamma_i$, $1\leq
i\leq n$ in a decomposition of $\bar F$ into projective planes are
disjoint from $E$, as are all the Dehn twists and cross-cap slides
we perform on $\bar F$ in the proof of Theorem~\ref{main}. Hence
the Dehn twists and cross-cap slides we perform give
homeomorphisms of $F$ which are the identity on the boundary
components.

Let $\alpha_i=[\gamma_i]$ and let $\bar\alpha_i$ be the images of
these elements in $H_1(\bar F,\Z)$. By choosing appropriate
orientations, we get that $H_1(F,\Z)$ is generated by the elements
$\alpha_i$ and $\beta_j$ with the relation
\begin{equation}
2\sum_i \alpha_i=\sum_j\beta_j.
\end{equation}
Note that as $H_1(\bar F,\Z)=H_1(F,\Z)/{\langle \beta_j\rangle}$,
it follows by the hypothesis that $\varphi$ induces an automorpism
$\bar\varphi$ of $H_1(\bar F,\Z)$. By Theorem~\ref{main} (and its
proof), this is induced by a composition of Dehn twists and
cross-cap slides, hence a homeomorphism $g\hbox{\rm :}\ F\rightarrow F$.
By composing $\varphi$ by $g_*^{-1}$, we can assume that
$\bar\varphi$ is the identity.

Similarly, we can use homeomorphisms supported in $E$ (which do
not change any $\alpha_i$) to reduce to the case when the
permutation $\sigma$ is the identity, i.e.
$\varphi(\beta_j)=\epsilon_j\beta_j$. As
$\bar\varphi(\bar\alpha_j)=\bar\alpha_j$, we get
$\varphi(\alpha_i)=\alpha_i+\sum_j c_{ij}\beta_j$ for some
integers $c_{ij}$. We define the complexity of $\varphi$ to be
$C(\varphi)=\sum_{i,j}|c_{ij}|$.

If $\varphi$ is not the identity, we shall reduce the complexity
of $\varphi$ using homeomorphisms called {\it boundary
slides}~\cite{Ko} similar to cross-cap slides.
\end{proof}

\begin{lem}
There are homeomorphisms $h_{ij}$ of $F$ such that the induced
automorphism of $H_1(F,\Z)$ takes $\alpha_i$ to
$\alpha_i-\beta_j${\rm ,} maps $\beta_j$ to $-\beta_j$ and fixes
all other $\alpha$'s and $\beta$'s.
\end{lem}

\begin{proof}
We shall use boundary slides~\cite{Ko} of the surface $F$. Namely,
suppose $\alpha$ is an orientation reversing simple closed curve
on a surface $S'$ and $D$ is a small disc centered around a point
on $\alpha$. Let $S$ be the surface obtained by deleting the
interior of $D$. Consider a homeomorphism of $S'$ which is the
identity outside a neighbourhood of $\alpha$ and which is obtained
by dragging $D$ once around $\alpha$ so that $D$ is mapped to
itself. By construction this extends to a homeomorphism of $S$,
which we call a boundary slide.

As in the case of cross-cap slides, the automorphism of
$H_1(F,\Z)$ induced by the boundary slide of the boundary
component corresponding to $\beta_j$ along the simple closed curve
$\gamma_i$ (representing $\alpha_i$) is as in the statement of the
lemma.\hfill $\Box$
\end{proof}

Now suppose $\varphi$ is not the identity. Observe that as
$\varphi$ is a homomorphism, $2\sum_i
\varphi(\alpha_i)=\sum_j\varphi(\beta_j)$. Using
$\varphi(\alpha_i)=\alpha_i+\sum_j c_{ij}\beta_j$,
$\varphi(\beta_j)=\epsilon_j\beta_j$ and $2\sum_i
\alpha_i=\sum_j\beta_j$, we see that $\sum_j
c_{ij}\beta_j=(\epsilon_j-1)\beta_j$. As the elements $\beta_j$,
$1\leq j\leq n$ are independent, it follows that for each $j$,
$\sum_i c_{ij}=\epsilon_j-1$.

We now consider two cases. Firstly, if some $\epsilon_j=-1$, then
observe that postcomposing with $h_{ij}$ takes $\varphi(\alpha_i)$
to $\varphi(\alpha_i)-\varphi(\beta_j)=\varphi(\alpha_i)+\beta_j$.
Hence $c_{ij}$ is changed to $c_{ij}+1$ (and no other $c_{kl}$ is
changed). In particular, if $c_{ij}<0$, the complexity is reduced.
But as $\sum_i c_{ij}=\epsilon_j-1=-2$, we must have some
$c_{ij}<0$, and hence a move reducing complexity.

Suppose now that each $\beta_j$ is $1$. Then as $\sum_i
c_{ij}=\epsilon_j-1=0$, either each $c_{ij}=0$, in which case we
are done, or some $c_{ij}>0$. Observe that postcomposing with
$h_{ij}$ takes $\varphi(\alpha_i)$ to
$\varphi(\alpha_i)-\varphi(\beta_j)=\varphi(\alpha_i)-\beta_j$.
Hence $c_{ij}$ is changed to $c_{ij}-1$ (and no other $c_{kl}$ is
changed), and hence the complexity is reduced. Thus in finitely
many steps, we reduce to the case where $\varphi$ is the identity.\hfill $\Box$

\section*{Acknowledgements}

We would like to thank Shreedhar Inamdar for helpful
conversations.

\end{document}